\documentclass[a4paper,10pt]{amsart}
\usepackage{times}
\usepackage{amsmath,amsthm,amssymb,enumerate}
\usepackage{epsfig}
\usepackage{amssymb}
\usepackage{amsmath}
\usepackage{amssymb}
\usepackage{amsmath,amsthm}
\usepackage[latin1]{inputenc}
\usepackage{bm}
\usepackage{bbm}
\usepackage{esint}
\usepackage{tikz-cd}
\usetikzlibrary{decorations.pathreplacing}
\usepackage{amsfonts}
\usepackage{amsxtra}
\usepackage{euscript,mathrsfs}
\usepackage{color}
\usepackage[left=3.1cm,right=3.1cm,top=4cm,bottom=3cm]{geometry}
\usepackage[citecolor=blue,colorlinks=true]{hyperref}
\allowdisplaybreaks

\usepackage{tikz-cd}
\usetikzlibrary{decorations.pathreplacing}

\usepackage{enumitem}
\setenumerate{label={\rm (\alph{*})}}

\usepackage{amsfonts}
\usepackage{amsxtra}

\newtheorem{theorem}{Theorem}

\newtheorem{proposition}[theorem]{Proposition}

\usepackage{xcolor}
\colorlet{ColorPink}{red!30}
\usepackage{graphicx}

\newcommand{\R}{\mathbb R}
\newcommand{\Rn}{{\mathbb{R}}^{n}}
\newcommand{\RN}{\mathbb{R}^{N}}
\newcommand{\WW}{\mathrm{W}}
\newcommand{\Leb}{\mathscr{L}^{n}}
\newcommand{\wstar}{\stackrel{\ast}{\rightharpoonup}}

\newcommand{\F}{\mathscr{F}}

\newcommand{\dd}{\mathrm{d}}

\newcommand{\CC}{\mathrm{C}}
\newcommand{\restrict}{\begin{picture}(10,8)\put(2,0){\line(0,1){7}}\put(1.8,0){\line(1,0){7}}\end{picture}}

\DeclareMathOperator{\BV}{BV}

\newcommand{\mm}{\mathbb{M}}

\newcommand{\sobo}{\operatorname{W}}

\newcommand{\bv}{\operatorname{BV}}

\newcommand{\dashint}{\fint}

\allowdisplaybreaks

\allowdisplaybreaks

\begin{document}


\title[Regularity for Quasiconvex Problems]{Regularity for higher order quasiconvex problems \\ with linear growth from below}
\author[F.~Gmeineder \& J.~Kristensen]{Franz Gmeineder and Jan Kristensen}
\address[F.~Gmeineder]{Mathematical Institute, University of Bonn, Endenicher Allee 60, 53116 Bonn, Germany}
\address[J.~Kristensen]{Mathematical Institute, University of Oxford, Andrew Wiles Building, OX2 6GG, Oxford, United Kingdom}

\begin{abstract}
We announce new existence and $\varepsilon$-regularity results for minimisers of the relaxation of strongly quasiconvex integrals that on
smooth maps $u \colon \Omega \subset \Rn \to \RN$ are defined by
$$
u \mapsto \int_{\Omega} \! F(\nabla^{k}u(x)) \, \dd x .
$$
The results cover the case of integrands $F$ with $(1,q)$ growth in the full range of exponents $1<q<\tfrac{n}{n-1}$ for which a
measure representation of the relaxed functional is possible and the minimizers belong to the space $\BV^{k}$ of maps whose $k$-th
order derivatives are measures. 
\end{abstract}

\maketitle


A key theme in the vectorial calculus of variations is the identification of necessary and sufficient conditions that ensure the existence and (partial)
regularity of minima of energy functionals. The challenges coming from applications, notably mathematical elasticity theory, mean that one must
consider a broader variety of integrands (energy densities) than are normally considered in the standard theory. In this connection various basic
problems emerge, including how the functionals should be defined on the natural spaces where a minimizer is sought. Here we shall adopt the approach
taken by \textsc{Marcellini} \cite{marcellini} for quasiconvex integrals in the Sobolev context and define the functionals by relaxation. For further
applications and results for variational problems under nonstandard growth conditions we refer to \textsc{Mingione} \cite{GM} and \textsc{Zhikov}
\cite{zhikov} and the references therein.

If we let $\Omega$ be a bounded and open smooth domain in $\Rn$ and $k\in\mathbb{N}$, such functionals are, ignoring
for ease of exposition lower order terms and other less relevant dependencies, typically given on smooth maps $u \colon \Omega \to \RN$ by
\begin{align}\label{eq1}
\mathrm{I}[u , \Omega]:=\int_{\Omega} \! F(\nabla^{k}u)\, \dd x
\end{align}
where $F \colon \mm \to \R$ is a continuous integrand defined on the space $\mm := \odot^{k}(\Rn, \RN )$ of symmetric $k$-linear $\R^{N}$-valued maps on $\R^{n}$.
The existence of minima subject to Dirichlet boundary data, that we recall in the $k$-th order context means zero order trace on $\partial \Omega$ together with normal
derivatives on $\partial \Omega$ up to and including order $k-1$, is then achieved using the direct method by virtue of growth and quasiconvexity properties of $F$.
Assuming $F$ satisfies such conditions in a suitably strict form, it is then natural to inquire as to which regularity properties
minimizers share. In this note we announce some new results in this program, and proceed by discussing our set-up in detail.
\bigskip

\noindent
\textbf{The variational problem and existence of minimizers.}
The space $\mm := \odot^{k}(\Rn, \RN )$ of symmetric $k$-linear $\R^{N}$-valued maps on $\R^{n}$ is equipped with the usual operator norm corresponding to
the standard euclidean norms on the underlying spaces $\Rn$ and $\RN$. All these norms are denoted by $| \cdot |$, the precise meanings being clear from context.
We consider a continuous integrand $F \colon \mm \to \R$ satisfying the \textit{$(1,q)$ growth condition}:
\begin{equation}\label{(1,q)}
c_{1}|z| -c_{2} \leq F(z) \leq c_{3} \bigl( |z|+1 \bigr)^{q} \quad \forall z \in \mm ,
\end{equation}
where $c_i >0$ are positive constants and $q \geq 1$ a fixed exponent. It is well-known that quasiconvexity in the sense of Morrey \cite{Morrey}, and its higher
order variant considered first by Meyers \cite{Meyers}, together with suitable growth conditions, are closely related to lower semicontinuity and coercivity
(see \cite{cycK,GK1} for the latter in the case when the lower bound in (\ref{(1,q)}) is omitted). Here we recall that $F$ is said to be \textit{quasiconvex} if
for all $z \in \mm$ and all smooth compactly supported maps $\varphi \colon \Rn \to \RN$ we have
\begin{equation}\label{qc}
\int_{\Rn} \bigl( F(\nabla^{k} \varphi (x)+z) - F(z) \bigr) \, \dd x \geq 0.
\end{equation}  
Let $g \in \WW^{k,q}( \Rn , \RN )$ and consider the problem of minimizing the variational integral $\mathrm{I}[u, \Omega ]$ defined at (\ref{eq1})
over $u \in \WW^{k,q}_{g}( \Omega , \RN ) := g+ \WW^{k,q}_{0}(\Omega , \RN )$. Clearly any minimizing sequence $( u_j ) \subset \WW^{k,q}_{g}$ for $\mathrm{I}[ \cdot , \Omega ]$
will in view of the lower bound in (\ref{(1,q)}) and Poincar\'{e}'s inequality be bounded in $\WW^{k,1}$. It is well-known that this forces us to consider the
minimization over the larger space $\BV^k = \BV^{k}( \Omega , \RN )$ of maps whose derivatives of order $k$ are bounded measures.
(See for instance \cite{Dorronsoro} for basic properties of such maps.) We then find a subsequence (not relabelled) and a map $u \in \BV^{k}$
so that $u_{j} \wstar u$ in $\BV^{k}$, meaning that $u_{j} \to u$ strongly in $\WW^{k-1,1}$ and $\nabla^{k}u_j \Leb \restrict \Omega \wstar D^{k}u$ in the sense of bounded
measures on $\Omega$. Here $D^{k}u$ is the full distributional $k$-th derivative of $u$, which is a bounded $\mm$-valued Radon measure on $\Omega$. Its
Lebesgue-Radon-Nikod\'{y}m decomposition with respect to $\Leb$ is denoted by $D^{k}u = \nabla^{k}u \Leb \restrict \Omega + D^{k}_{s}u$.
An additional complication here is that the trace map is not continuous in the weak$\mbox{}^{\ast}$ topology
of $\BV^k$ and so, following \cite{GiMoSo}, it is useful to incorporate the Dirichlet boundary condition in the variational problem as a kind of penalization term.
For exponents $q>1$ this procedure leads to a rather implicit functional and penalization term. To define it we fix a bounded open domain $\Omega^{\prime}$ in $\Rn$
with $\Omega \Subset \Omega^{\prime}$. For each $u \in \BV^{k} (\Omega^{\prime} , \RN )$ define 
$$
\F_{g} [ u ] = \left\{
\begin{array}{ll}
\displaystyle{\int_{\Omega^{\prime}} \! F(\nabla^{k}u) \, \dd x} & \mbox{ when } u \in \WW^{k,q}(\Omega^{\prime},\RN ) \mbox{ and }
                                                                        u = g \mbox{ on } \Omega^{\prime}\setminus \Omega\\
\infty   & \mbox{ otherwise.}
\end{array}
\right.
$$
The functional to be minimized is now the \textit{sequential weak$\mbox{}^{\ast}$ lower semicontinuous envelope} of $\F_{g}[ \cdot ]$ on $\BV^{k}( \Omega^{\prime},\RN )$:
$$
\overline{\F}_{g}[ u ] = \inf\left\{ \liminf_{j\to\infty} \F_{g} [u_{j}] \colon
\begin{array}{c}
(u_{j})\subset \BV^{k}( \Omega^{\prime},\RN ),\\
u_{j}\wstar u\;\text{in}\;\BV^{k}( \Omega^{\prime},\RN ) 
\end{array}\right\}
$$
It is evident that $\overline{\F}_{g}[ \cdot ]$ admits minimizers on $\BV^{k}( \Omega^{\prime},\RN )$ and these minimizers will all agree with $g$
on $\Omega^{\prime}\setminus \Omega$. For $q=1$ we have a formula for $\overline{\F}_{g}[ \cdot ]$ that is a straightforward generalization of the formula
that is known in the first order case $k=1$ (see \cite{AD,FoMu,KR}) and we shall therefore in the following mainly focus on exponents $q>1$ where
there is no known formula for the relaxed functional. Adapting arguments from \cite{FMParma} we extend the results of \cite{JK1} and
\cite{schmidt3} and show the following:

\begin{proposition}\label{relaxedF}
Assume $F \colon \mm \to \R$ is continuous and satisfies (\ref{(1,q)}), (\ref{qc}) for an exponent $q \in [1,\tfrac{n}{n-1})$.
For $u \in \WW^{k,q}(\Omega^{\prime},\RN )$ with $u = g$ on $\Omega^{\prime}\setminus \Omega$ we have
$$
\overline{\F}_{g}[ u ] = \int_{\Omega^{\prime}} \! F( \nabla^{k}u) \, \dd x.
$$
Furthermore, if $\bar{u} \in \BV^{k}( \Omega^{\prime} ,\RN )$ is a minimizer for $\overline{\F}_{g}[ \cdot ]$, then for balls $B=B_{r}(x_{0}) \Subset \Omega$
we have
$$
\int_{B} \! F(\nabla^{k}\bar{u}) \, \dd x \leq \int_{B} \! F(\nabla^{k}\bar{u}+ \nabla^{k}\varphi ) \, \dd x
$$
for all $\varphi \in \CC^{k}_{c}(B,\RN )$. 
\end{proposition}
When the integrand $F$ is differentiable, the last result of the proposition together with standard arguments yield an Euler-Lagrange equation for
the minimizer $\bar{u}$. We emphasize that here $\nabla^{k}\bar{u}$ is merely the $\Leb$ density of the full distributional derivative $D^{k}\bar{u}$
and as such does not have gradient structure in general. The arguments rely on extending the measure representation results from \cite{FMParma} in the $(p,q)$ growth case
to the $(1,q)$ growth case and for this the restriction $q < \tfrac{n}{n-1}$ is essential as demonstrated in \cite{AcerbiDalMaso}. We also remark
that the first result of the proposition, that the relaxed functional $\overline{\F}_{g}[ \cdot ]$ is an extension of $\F_g$ from $\WW^{k,q}$ maps agreeing
with $g$ off $\Omega$, fails in a related situation of $(p,n)$ growth when $p<n-1$, see \cite{maly}. We finally mention that the \textit{pointwise definition}
of the variational integrals that is used in \cite{BM}, in the $k$-th order would require the stronger condition of $\WW^{k,1}$ quasiconvexity
and one would still need to make a relaxation from $\WW^{k,1}$ to $\BV^{k}$ that incorporates the Dirichlet boundary condition.
However, in this case the relaxation formula is known in the related first order case, $k=1$, and integrands of $(p,q)$ growth with $1<p \leq q < \tfrac{np}{n-1}$
(see \cite{JK2}).
\bigskip

\noindent
\textbf{Regularity of minimizers}
The main result of \cite{GK2} concerns the regularity of minimizers of the relaxed functional $\overline{\F}_{g}[ \cdot ]$ defined above under an additional
natural strong quasiconvexity assumption on the integrand. To define it let
$$
E(z) = \sqrt{1+|z|^{2}}-1, \quad z \in \mm
$$
be the \textit{reduced area integrand} as also considered in \cite{GK1}. Then we say that
the continuous integrand $F \colon \mm \to \R$ is \textit{strongly quasiconvex} if there exists a positive constant $\ell >0$ such that $z \mapsto F(z)-\ell E(z)$
is quasiconvex:
\begin{equation}\label{sqc}
\int_{\Rn} \bigl( F(z+\nabla^{k}\varphi (x) ) -F(z) \bigr) \, \dd x \geq \ell\int_{\Rn} \bigl( E(z+\nabla^{k}\varphi (x) ) -E(z) \bigr) \, \dd x
\end{equation}
holds for all $z \in \mm$ and all $\varphi \in \CC^{\infty}_{c}(\Rn , \RN )$. It is known that in the $p$-growth case higher order quasiconvexity is related
to (first order) quasiconvexity (see \cite{sverak}, \cite{DFLM} and \cite{Cagnetti}), and we show in \cite{GK3} that these results persist in the linear
growth case.

When $u$ is of class $\BV^{k}$ we write for a ball $B=B_{r}(x_{0})$, 
$$
(Du)_{B} = \frac{Du (B)}{\Leb (B)}
$$
and
$$
\dashint_{B} \! E(D^{k}u - (D^{k}u)_{B}) := \frac{1}{\Leb (B)}\left( \int_{B} \! E(\nabla^{k}u(x)-(D^{k}u)_{B})
\, \dd x + |D^{k}_{s}u|(B) \right) ,
$$
where $|D^{k}_{s}u|$ is the total variation measure for the singular $\mm$-valued measure $D^{k}_{s}u$. In these terms we have

\begin{theorem}\label{thm:main1}
Let $\overline{\F}_{g}[ \cdot ]$ be the relaxed functional defined above and assume the integrand $F \colon \mm \to \R$ satisfies
\begin{itemize}
\item[(H1)] $F$ is $\CC^{2,1}_{\mathrm{loc}}(\mm )$,
\item[(H2)] $F$ has $(1,q)$ growth for some exponent $q \in [1, \tfrac{n}{n-1})$,
\item[(H3)] $F$ is strongly quasiconvex.
\end{itemize}
Then for each $m>0$ there exists $\varepsilon_{m}>0$ depending on the data in (H1), (H2), (H3) with the following property. For any
minimizer $\bar{u} \in \BV^{k}( \Omega^{\prime}, \RN )$ of $\overline{\F}_{g}[ \cdot ]$ and each ball $B=B_{r}(x_{0}) \subset \Omega$
with
$$
|(D\bar{u})_{B}| < m \quad \mbox{ and } \quad \dashint_{B} \! E(D\bar{u}-(D\bar{u})_{B}) < \varepsilon_{m}
$$
we have that $\bar{u}$ is $\CC^{k+1,\alpha}_{\mathrm{loc}}$ on $B_{r/2}(x_{0})$ for all $\alpha <1$. In particular, $\bar{u}$ is partially
$\CC^{k+1,\alpha}$ regular in $\Omega$ for each $\alpha <1$: there exists a relatively closed subset $\Sigma_{\bar{u}} \subset\Omega$ with $\Leb (\Sigma_{\bar{u}})=0$
such that $\bar{u}|_{\Omega \setminus \Sigma_{\bar{u}}}$ is locally $\CC^{k+1,\alpha}$ for each $\alpha < 1$.
\end{theorem}
The case $k=q=1$ was proved, in a more general set-up, in \cite{GK1}. We also remark that the result is new also when we strengthen the
hypothesis (H3) to \textit{strong convexity} meaning that $z \mapsto F(z)-\ell E(z)$ is convex for some $\ell > 0$. We also remark that, under
slightly more restrictive bounds on the exponent $q$ (as in \cite{ABF}), we can in the strongly convex case allow the integrand to depend
on $x \in \Omega^{\prime}$, whereby several examples exhibiting $(p,q)$ growth with $p>1$ treated in the literature are extended to the
$(1,q)$ growth scenario.

The results of \textsc{Schmidt} \cite{schmidt3} (see also \cite{schmidt1,schmidt2}) concern the first order case and integrands
of $(p,q)$ growth with $p>1$, and more restrictive bounds on the exponent $q$, but also include certain degenerate situations not covered here.
These results were preceeded by work of \textsc{Fusco \& Hutchinson} \cite{FuHu} on the anisotropic polyconvex case, and certain related
quasiconvex situations (see also \cite{EsMi}).
Finally \textsc{Schemm} \cite{schemm} extended the results of \cite{schmidt1} to the $k$-th order case. We refer to \cite{GK1} for a
discussion of further regularity results for $\BV$ minimizers of convex and quasiconvex integrals, and to \cite{GM} for a survey
of regularity results in Calculus of Variations and the related PDEs in the Sobolev context.

\end{document}